\documentclass[12pt]{amsart}
\usepackage{amscd}
\def\NZQ{\Bbb}               
\def\NN{{\NZQ N}}

\def\ZZ{{\NZQ Z}}

\def\D{{\Delta}}
\def\mm{{\frak m}}
\def\F{{\mathcal F}}
\def\a{{\bold a}}
\def\b{{\bold b}}
\def\c{{\bold c}}

\def\e{{\bold e}}
\def\1{{\mathbf 1}}
\def\0{{\mathbf 0}}
\def\opn#1#2{\def#1{\operatorname{#2}}} 
\opn\diam{diam} 
\opn\height{ht}

\newtheorem{Theorem}{Theorem}[section]
\newtheorem{Lemma}[Theorem]{Lemma}
\newtheorem{Corollary}[Theorem]{Corollary}

\newtheorem{Example}[Theorem]{Example}

\begin{document}

\title{Cohen-Macaulayness of powers\\ of two-dimensional squarefree monomial ideals}

\author{Nguyen Cong Minh}
\address{Department of Mathematics, University of Education, 136 Xuan Thuy, Hanoi, 
Vietnam}
\email{ngcminh@@gmail.com}

\author{Ngo Viet Trung}
\address{Institute of Mathematics, 18 Hoang Quoc Viet, Hanoi, 
Vietnam}
\email{nvtrung@@math.ac.vn}

\begin{abstract}
Two-dimensional squarefree monomial ideals can be seen as the Stanley-Reisner ideals of graphs.
The main results of this paper are combinatorial characterizations for the Cohen-Macaulayness of ordinary and symbolic powers of such an ideal in terms of the associated graph. 
\end{abstract}
\maketitle

\section*{Introduction}

Let $I$ be a two-dimensional squarefree monomial ideal in a polynomial ring $k[X] = k[x_1,...,x_n]$, where $k$ is a field and $n \ge 3$.
Since we are studying the Cohen-Macaulayness of powers of $I$, we will assume that $I$ is unmixed.
In this case, $I$ is the intersection of prime ideals of the form
$P_{ij} := \big(X \setminus \{x_i,x_j\}\big).$
If we think of the indices of such ideals as the edges of a graph $G$ on the vertex set $\{1,...,n\}$, then
$$I = \bigcap_{\{i,j\} \in G}P_{ij}.$$
In other words, $I$ is the Stanley-Reisner ideal of the one-dimensional simplicial complex associated with $G$  \cite{BH}, \cite{Sta}.
\smallskip

It is well-known that $I$ is a Cohen-Macaulay ideal if and only if $G$ is a connected graph. 
Our original goal is to characterize the Cohen-Macaulayness of $I^m$ in terms of $G$ for each $m \ge 2$.
This kind of problem is usually hard and there are few classes of ideals for which one knows exactly which powers of them are Cohen-Macaulay. Notice that $I^m$ is Cohen-Macaulay for all $m \ge 2$ if and only if $I$ is a complete intersection \cite{AV}, \cite{Wa} (A similar characterization for the generalized Cohen-Macaulayness of $I^m$ can be found in \cite{GT}). Since the Cohen-Macaulayness of $I^m$ implies the equality $I^{(m)} = I^m$, where $I^{(m)}$ denotes the $m$-th symbolic power of $I$, we have to study the problems when $I^{(m)}$ is a Cohen-Macaulay ideal and when $I^{(m)} = I^m$. These problems are of independent interest and they will be solved completely in this paper.
\smallskip

For the Cohen-Macaulayness of $I^{(m)}$ we need to distinguish two cases $m = 2$ and $m \ge 3$.
We shall see that $I^{(2)}$ is Cohen-Macaulay if and only if $\diam(G) \le 2$,
where $\diam(G)$ denotes the diameter of $G$. For $m \ge 3$ we show that 
$I^{(m)}$ is Cohen-Macaulay if and only if every pair of disjoint edges of $G$ is contained in a cycle of length 4. In that case, all other symbolic powers of $I$ are Cohen-Macaulay, too. The proof is based on a formula for local cohomology modules of quotient rings of monomial ideals found by Takayama \cite{Ta}. 
\smallskip

The equality between ordinary and symbolic powers is a purely ideal-theoretic problem. In turns out that there are few graphs which satisfies the condition $I^{(m)} = I^m$ for some $m \ge 2$. In fact, we shall see that $I^{(2)} = I^2$ if and only if $n =3,4,5$ and $G$ is a path or a cycle or the union of two disjoint edges, and for $m \ge 3$, $I^{(m)} = I^m$  if and only if $n = 3,4$ and $G$ is a path or a cycle or the union of two disjoint edges.
\smallskip

As a consequence of the above results, we show that $I^2$ is Cohen-Macaulay if and only if $n = 3$ or $G$ is a cycle of length $4,5$, and for $m \ge 3$, $I^m$ is Cohen-Macaulay  if and only if $n = 3$ or $G$ is a cycle of length 4. In particular, $I$ must be a complete intersection if $I^m$ is Cohen-Macaulay for some $m \ge 3$.
\smallskip 

The results on the equality between ordinary and symbolic powers can be also used to study vertex cover algebras of simplicial complexes recently introduced in \cite{HHT}. We shall see that the vertex cover algebra of a pure $(n-3)$-dimensional simplicial complex $\D$ on $n$ vertices is standard graded if and only if $n = 3$ or $n = 4$ and $\D$ is a path or a cycle or the union of two disjoint edges. 

\section{Takayama's formula}

Let $S = k[X] = k[x_1,...,x_n]$. Let $I$ be an arbitrary monomial ideal in $S$. Since $S/I$ is an $\NN^n$-graded algebra,
$H_\mm^i(S/I)$ is a $\ZZ^n$-graded module over $S/I$.
For every degree $\a \in \ZZ^n$ we denote by $H_\mm^i(S/I)_\a$ the $\a$-component of $H_\mm^i(S/I)$.
Inspired of a result of Hochster in the squarefree case \cite[Theorem 4.1]{Ho}, Takayama found the following combinatorial formula for $\dim_kH_\mm^i(S/I)_\a$ \cite[Theorem 2.2]{Ta}.
\smallskip

For every $\a = (a_1,...,a_n) \in \ZZ^n$ we set $x^\a = x_1^{a_1}\cdots x_n^{a_n}$ and $G_\a = \{i|\ a_i <0\}$. We denote by $\D_\a(I)$ the simplicial complex of all $F \subseteq \{1,...,n\}$ such that \par
(1) $F \cap G_\a = \emptyset$, \par
(2) for every minimal generator $x^\b$ of $I$ there exists an index $i \not\in F \cup G_\a$ with $b_i > a_i$.
\smallskip

Let $\D(I)$ denote the simplicial complex of all $F$ such that $\prod_{i \in F}x_i\not\in \sqrt{I}$.

\begin{Theorem}[Takayama's formula]\label{Takayama}
$$\dim_kH_\mm^i(R/I)_\a = 
\begin{cases}
\dim_k\tilde H_{i-|G_\a|-1}(\D_\a(I),k) & \text{\rm if }\ G_\a \in \D(I),\\
0 & \text{\rm else. }
\end{cases} $$
\end{Theorem}
 
Notice that the original formula in \cite{Ta} is a bit different. 
It puts beside $G_\a \in \D(I)$ the condition $a_i \le \rho_i -1$ for $i = 1,...,n$, where $\rho_i$ is the maximum of the $i$th coordinates of all vectors $\b \in \NN^n$ such that $x^\b$ is a minimal generator of $I$. But the proof in \cite{Ta} shows that we may drop this condition, which is more convenient for our investigation. 
\smallskip

Using Takayama's formula we obtain the following condition for the vanishing of the first local cohomology module.

\begin{Corollary}\label{H1}
$H_\mm^1(S/I) = 0$ if and only if $\D_\a(I)$ is connected for all $\a \in \NN^n$.
\end{Corollary}

\begin{proof}
We need to check when $\tilde H_{-|G_\a|}(\D_\a(I),k) = 0$.
Therefore, we only need to consider the case $|G_\a| = 0$ or, equivalently, $\a \in \NN^n$.
The conclusion comes from the fact that $\tilde H_0(\D_\a(I),k) = 0$ if and only if $\D_\a(I)$ is connected.
\end{proof}

As a consequence, if we want to study the Cohen-Macaulayness of a two-dimensional monomial ideal $I$, we only need to check the connectedness of $\D_\a(I)$ for $\a \in \NN^n$.  We shall see below that there are only a finite number of possibilities for $\D_\a(I)$. For that we need to reformulate the definition of $\D_\a(I)$ in a more simple way. 
\smallskip

For every subset $F$ of $\{1,...,n\}$ let $S_F = S[x_i^{-1}|\ i \in F \cup G_\a]$ and 
$I_F =IS_F$. It is easy to check that $x^\a \in S_F$ and $\D_\a(I)$ is the simplicial complex of all $F$ such that $F \cap G_\a = \emptyset$ and $x^\a \not\in I_F$. If $\a \in \NN^n$, we may replace $I_F$ by the ideal obtained from $I$ by setting $x_i = 1$ for $i \in F$.  
\smallskip

Using this interpretation of $\D_\a(I)$ we can show that $\D_\a(I)$ is closely related to $\D(I)$.
Notice first that if we denote by $P_F$ the ideal generated by the variables $x_i$, $i \not\in F$, then $$\sqrt{I}= \bigcap_{F \in \F(\D(I))}P_F,$$
where $\F(\D(I))$ denotes the set of the facets of $\D(I)$. 

\begin{Lemma}\label{subcomplex}
$\D_\a(I)$ is a subcomplex of $\D(I)$. Moreover, if $I$ has no embedded associated prime ideals and $\a  \in \NN^n$, the facets of $\D_\a(I)$ are facets of $\D(I)$.
\end{Lemma}

\begin{proof} 
Let $F$ be an arbitrary face of $\D_\a(I)$. If $F \not\in \D(I)$, then $\prod_{i\in F}x_i \in \sqrt{I}$. Therefore, $\sqrt{I_F} = S_F$, which contradicts the condition $x^\a \not\in I_F$. 
For the second assertion we only need to show that $F$ is contained in a facet of $\D(I)$ which belongs to $\D_\a(I)$.
Observe first that $I_F = JS_F$, where $J$ is the intersection of all primary components of $I$ whose associated primes do not contain any variable $x_i$ with $i \in F$. Since $x^\a \not\in I_F$, there is at least a primary component $Q$ of $I$ such that $x^\a \not\in Q$. Let $G$ be the set of indices such that $P_G$ is the associated prime of $Q$. Then $F \subseteq G$. If $I$ has no embedded associated prime ideals, then $P_G$ is a prime component of $\sqrt{I}$, hence $G$ must be a facet of $\D(I)$. If $\a \in \NN^n$, then $G_\a = \emptyset$. Since $x^\a \not\in QS_G = I_G$, we can conclude that $G \in \D_\a(I)$. 
\end{proof}

\begin{Example}\label{example 1} {\rm \noindent (1) For $\a = 0$ we have $\D_\a(I) = \D(I)$ because for all facets $F$ of $\D(I)$ we have $I_F \neq S_F$ and therefore $x^\a = 1 \not\in I_F$ for all facets $F$ of $\D(I)$.\par
\noindent (2) For $\a = \e_i$, where $\e_i$ denotes the $i$th unit vector of $\NN^n$, then the facets of $\D_\a(I)$ are 
the facets $F$ of $\D(I)$ with the property that $x_i$ is not contained in the primary component of $I$ associated with $P_F$. In particular, if no primary components of $I$ contains $x_i$, then $\D_\a(I) = \D(I)$.}
\end{Example}

The above lemma will play a crucial role in the study of the symbolic powers of two-dimensional squarefree monomial ideals in the next section.

\section{Cohen-Macaulayness of symbolic powers}

Let $G$ be a simple graph (i.e. without isolated vertices and loops) on the vertex set $\{1,...,n\}$, $n \ge 3$.
We associate with $G$ the ideal
$$I_G := \bigcap_{\{i,j\} \in G} P_{ij},$$
where $P_{ij} = (x_t|\ t \neq i,j)$. By definition, the $m$-th symbolic power of $I_G$ is the ideal
$$I_G^{(m)} := \bigcap_{\{i,j\} \in G} P_{ij}^m.$$

If $n = 3$, $I_G$ is a principal ideal. Hence, $I_G^{(m)} = I_G^m$  and $I_G^m$ is Cohen-Macaulay  for all $m \ge 2$. For $n \ge 3$, we will use Takayama's lemma to study the Cohen-Macaulayness of $I_G^{(m)}$.
\smallskip

It is obvious that $\D(I_G^{(m)}) = \D(I_G)$ for all $m \ge 2$ and that $\D(I_G)$ is the one-dimensional simplicial complex associated with $G$. For simplicity, we identify $\D(I_G)$ with $G$ so that the facets of $\D(I_G)$ are edges of $G$. By Lemma \ref{subcomplex}, $\D_\a(I_G^{(m)})$ is a simple subgraph of $G$ for all $\a \in \NN^n$. This graph can be described as follows.

\begin{Lemma}\label{Delta}
$\D_\a(I_G^{(m)})$ is the simple graph of the edges 
$\{i,j\} \in G$ with $\sum_{t \neq i,j}a_t < m$.
\end{Lemma}

\begin{proof}
We have to determine which edge $\{i,j\} \in G$ gives a facet of $\D_\a(I_G^{(m)})$. 
For $I = I_G^{(m)}$ and $F = \{i,j\}$ we have $I_F = P_{ij}^m$. Since $G_\a = \emptyset$ we only need to check when $x^\a \not\in P_{ij}^m$. But this condition is satisfied if and only if $\sum_{t \neq i,j}a_t < m$. 
\end{proof}

\begin{Example}
{\rm Let $m \ge 2$ and $\a = (m-1)(\e_r +\e_s)$, where $r < s$ are two arbitrary indices. For any edge $\{i,j\} \in G$ we have 
$$\sum_{t \neq i,j}a_t = 2(m-1) - a_i - a_j < m$$
if and only if the edge $\{i,j\}$ contains $r$ or $s$ or both. Therefore, $\D_\a(I_G^{(m)})$ is the subgraph of $G$ which consists of all edges containing $r$ or $s$ or both.}
\end{Example}

If $I_G^{(m)}$ is Cohen-Macaulay, then $\D_\a(I_G^{(2)})$ is connected for $\a = (m-1)(\e_r +\e_s)$ by Lemma \ref{H1}. By the description of $\D_\a(I_G^{(2)})$ in the above example, its connectivity implies that $\{r,s\} \in G$ or there exists an index $t$ such that $\{r,t\},\{s,t\} \in G$. In other words, the minimal length of paths from $r$ to $s$ is at most 2. Since $r,s$ can be any pair of vertices, $G$ must be connected.  
\smallskip

In graph theory, the distance between two vertices of $G$ is the minimal length of paths from one vertex to the other vertex. This length is infinite if there is no paths connecting them. The maximal distance between two vertices of $G$ is called the diameter  of $G$ and denoted by $\diam(G)$. 
\smallskip

By the above observation, $\diam(G) \le 2$ if $I_G^{(m)}$ is Cohen-Macaulay. It turns out that this is also a sufficient condition for the Cohen-Macaulayness of $I_G^{(2)}$.

\begin{Theorem} \label{CM 1}
$I_G^{(2)}$ is Cohen-Macaulay if and only if $\diam(G) \le 2$.
\end{Theorem}

\begin{proof}
As we have seen above, we only need to prove the sufficient part.
Assume that $\diam(G) \le 2$. If $I_G^{(2)}$ is not Cohen-Macaulay, there exists $a \in \NN^n$ such that $\D_\a(I_G^{(2)})$ is not connected by Lemma \ref{H1}.
Since $\D_\a(I_G^{(2)})$ is a simple subgraph of $G$, $\D_\a(I_G^{(2)})$ must contain two edges, say $\{1,2\}$ and $\{3,4\}$, such that their vertices are not connected by a path. By Lemma \ref{Delta} we have
$$\sum_{t \neq 1,2}a_t  < 2\;\; \text{\rm and}\;
\sum_{t \neq 3,4}a_t < 2.$$
The first inequality implies that among the components $a_t$ with $t \neq 1,2$ there is at most a component $a_t \neq 0$ and this component must be 1. Similarly, among the components $a_t$ with $t \neq 3,4$ there is at most a component $a_t \neq 0$ and this component must be 1. Therefore, we either have $\a = 0$ or $\a = \e_t$ for some $t$ or $\a = \e_r + \e_s$ with $r \in \{1,2\}$ and $s \in \{3,4\}$. \par
If $\a = 0$ or $\a = \e_t$, then $\D_\a(I_G^{(2)}) = G$ by Example \ref{example 1}, which is a contradiction because $G$ is connected.
 If $\a = \e_r + \e_s$, then $\D_\a(I_G^{(2)})$ is the subgraph of $G$ which consists of all edges containing $r$ or $s$ or both.
Since $\diam(G) \le 2$, either $\{r,s\} \in G$ or there exists an index $t$ such that
$\{r,t\},\{s,t\} \in G$. Therefore, $\D_\a(I_G^{(2)})$ is connected so that we get a contradiction, too.
\end{proof}

For the Cohen-Macaulayness of $I_G^{(m)}$, $m \ge 3$, we need a stronger condition on $G$.

\begin{Theorem}\label{CM 2}
For $m \ge 3$, $I_G^{(m)}$ is Cohen-Macaulay if and only if every pair of disjoint edges of $G$ is contained in a cycle of length 4. 
\end{Theorem}

\begin{proof}
For simplicity let $I = I_G^{(m)}$.
Assume that $I$ is Cohen-Macaulay and that $G$ has two disjoint edges, say $\{1,2\}$ and $\{3,4\}$.
Consider the vector $\a \in \NN^n$ with 
$$a_1 = 1, a_2 = 1, a_3 = m-1,\; \text{\rm and } a_i = 0\; \text{for } i = 4,...,n.$$
Then $\sum_{t \neq 1,2}a_t = m-1$ and $\sum_{t \neq 3,4}a_t = 2 < m$, hence $\{1,2\}, \{3,4\} \in \D_\a(I)$ by Lemma \ref{Delta}.
For $i = 4,...,n$ we have $\sum_{t \neq 1,i}a_t = m$ and $\sum_{t \neq 2,i}a_t = m$ so that $\{1,i\}, \{2,i\} \not\in \D_\a(I)$. 
Since $\D_\a(I)$ is a connected graph by Lemma \ref{H1}, there must be a path in $\D_\a(I)$ connecting the vertices $1$ and $3$. This is the case only when  $\{1,3\}$ or $\{2,3\}$ belongs to $\D_\a(I)$. From this it follows that $\{1,3\}$ or $\{2,3\}$ belongs to $G$. 
Similarly, we can show that every vertex of $\{1,2\}$ or $\{3,4\}$ is connected at least by an edge of $G$ with the other edge. Now it is easy to see that there are two other edges of $G$ with vertices in $\{1,2,3,4\}$ such that together with $\{1,2\}, \{3,4\}$ they form a cycle of length 4. \par
Conversely, assume that every pair of disjoint edges of $G$ is contained in a cycle of length 4. By Lemma \ref{H1}, if $I$ is not Cohen-Macaulay,  there exists $\a \in \NN^n$ such that $\D_\a(I)$ is not connected. Since $\D_\a(I)$ is a simple graph, $\D_\a(I)$ must contain two disjoint edges, say $\{1,2\}, \{3,4\}$ which belongs to different connected components of  $\D_\a(I)$. By Lemma \ref{Delta}, the edges $\{1,2\}, \{3,4\}$ belong to $G$ and
$$\sum_{t \neq 1,2}a_t < m,\; \sum_{t \neq 3,4}a_t < m .$$
We may assume that $\{1,3\}, \{2,4\} \in G$ in oder to get a cycle of length 4 containing the edges $\{1,2\}, \{3,4\}$. Since 
$\{1,3\}, \{2,4\} \not\in \D_\a(I)$, we have
$$\sum_{t \neq 1,3}a_t \ge m,\; \sum_{t \neq 2,4}a_t \ge m.$$
The above inequalities yields
$$2m > \sum_{t \neq 1,2}a_t + \sum_{t \neq 3,4}a_t = \sum_{t \neq 1,3}a_t + \sum_{t \neq 2,4}a_t \ge 2m,$$
a contradiction. So we can conclude that $I_G^{(m)}$ is Cohen-Macaulay.
\end{proof}

\begin{Corollary}
If $I_G^{(m)}$ is Cohen-Macaulay for some $m \ge 3$, then $I_G^{(m)}$ is Cohen-Macaulay for all $m \ge 1$.
\end{Corollary}

\begin{proof}
Since the condition on $G$ of Theorem \ref{CM 2} is independent of $m$, $I_G^{(m)}$ is Cohen-Macaulay for all $m \ge 3$ once $I_G^{(m)}$ is Cohen-Macaulay for some $m \ge 3$. Since this condition implies $G$ is connected with $\diam(G) \le 2$, it implies that $I_G$ and $I_G^{(2)}$ are Cohen-Macaulay, too.
\end{proof}

Using the above theorems we can easily construct a graph $G$ such that $I_G^{(2)}$ is Cohen-Macaulay but $I_G^{(m)}$ is not Cohen-Macaulay for all $m \ge 3$.

\begin{Example}
{\rm Let $G$ be a cycle of length 5. Then $\diam(G) = 2$ but no pair of edges of $G$ is contained in a cycle of length 4.}
\end{Example}

Notice that there are infinitely many graphs which satisfy the conditions of Theorem \ref{CM 1} and Theorem \ref{CM 2}. 
An instance is the class of complete graphs.

\section{Equality between ordinary and symbolic powers}

Let $G$ be a simple graph on the vertex set $\{1,...,n\}$, $n \ge 3$.
We want to study the equality between the ordinary and symbolic powers of the ideal
$$I_G := \bigcap_{\{i,j\} \in G} P_{ij},$$
where $P_{ij} = (x_t|\ t \neq i,j)$. Recall that the $m$-th symbolic power of $I_G$ is the ideal
$$I_G^{(m)} := \bigcap_{\{i,j\} \in G} P_{ij}^m.$$

We observed first that $I_G$ is generated by the monomials $x_ix_j$, where $\{i,j\}$ is a non-edge of $G$, and $x_ix_jx_k$, where $\{i,j,k\}$ is a triangle of $G$. 
\smallskip

If $n = 3$, then $G$ is either a path or a triangle. Hence $I_G$ is either of the form $(x_1x_2)$ or $(x_1x_2x_3)$, which satisfies the condition $I_G^{(m)} = I_G^m$ for all $m \ge 2$. 
\smallskip

If $n \ge 4$, there are only a few graphs which satisfies the condition $I_G^{(m)} = I_G^m$ for some $m \ge 2$. To see that we shall need the following observation.

\begin{Lemma} \label{triangle}
Let $n \ge 4$. If $G$ satisfies the condition $I_G^{(m)} = I_G^m$ for some $m \ge 2$, then $n = 4,5$ and $G$ is a path or a cycle or a union of two disjoint edges.
\end{Lemma}

\begin{proof}
We will first show that if $G$ has a triangle of edges or a triangle of non-edges, then $I_G^{(m)} \neq I_G^m$ for all $m \ge 2$.\par
If $G$ has a triangle of edges, say with the vertices $\{1,2,3\}$, then $x_1x_2, x_2x_3, x_1x_3 \not\in I_G$ and $x_1x_2x_3 \in I_G$. Since $x_1x_2x_3x_4 \in I_G^{(2)}$, we have 
$$x_1^{m-1}x_2^{m-1}x_3^{m-1}x_4 = (x_1x_2x_3)^{m-2}x_1x_2x_3x_4 \in I_G^{m-2}I_G^{(2)} \in I_G^{(m)}.$$
On the other hand, if $x_1^{m-1}x_2^{m-1}x_3^{m-1}x_4 \in I_G^m$, then it is divisible by a
product of $m$ monomials in $I_G$. Since only one of them contain $x_4$, the others $m-1$ monomials involve only $x_1,x_2,x_3$. Thus, the product of these $m-1$ monomials is divisible by $(x_1x_2x_3)^{m-1}$ so that the remained monomial must be $x_4$, a contradiction. \par
If $G$ has a triangle of non-edges, say with the vertices $\{1,2,3\}$, then $x_1x_2 \in I_G$ and $x_1x_2x_3 \in I_G^{(2)}$. Therefore, 
$$x_1^{m-1}x_2^{m-1}x_3 = (x_1x_2)^{m-2}(x_1x_2x_3) \in I_G^{m-2}I_G^{(2)} \subseteq I_G^{(m)}.$$
But $x_1^{m-1}x_2^{m-1}x_3 \not\in I_G^m$ because it has degree $2m-1$, whereas the minimal degree of the elements of $I_G^m$ is $2m$.\par

Therefore, if $G$ satisfies the condition $I_G^{(m)} = I_G^m$ for some $m \ge 2$, then $G$
has no triangles of edges and no triangles of non-edges. As a consequence, every vertex belongs to at most two edges and to at most two non-edges of $G$. For instance, if $\{1,2\},\{1,3\}, \{1,4\} \in G$, we must have  $\{2,3\},\{2,4\},\{3,4\} \not\in G$, which forms a triangle of non-edges, a contradiction.
Using these properties we can easily check that $n \le 5$ and that $G$ is a path or a cycle in the connected case or the union of two disjoint edges in the unconnected case. 
\end{proof}

It turns out that the necessary condition of Lemma \ref{triangle} is also a sufficient condition
in the case $m = 2$.  

\begin{Theorem} \label{equa 1}
Let $n \ge 4$. Then $I_G^{(2)} = I_G^2$ if and only if $n = 4,5$ and $G$ is a path or a cycle or the union of two disjoint edges.
\end{Theorem}

\begin{proof}
By Lemma \ref{triangle}, we only need to show the sufficient part. 
If $n = 4,5$ and $G$ is a path or a cycle or the union of two disjoint edges, then $G$ has no triangles of edges and no triangles of  non-edges. Using this condition we will show that every monomial $f\in I_G^{(2)}$ belongs to $I_G^2$. \par
If $f$ involves only two variables, say $x_1,x_2$, then every ideal $P_{ij}$ must contain $x_1$ or $x_2$.
Therefore, $\{1,2\} \not\in G$. Hence $x_1x_2 \in I_G$. Since $f \in P_{1j}^2$ for some $j \neq 2$ and since $P_{1j}$ does not contain $x_1$, $f$ is divisible by $x_2^2$. Similarly, $f$ is divisible by $x_1^2$. Hence $f$ is divisible by $(x_1x_2)^2$ so that $f \in I_G^2$.\par
If $f$ involves only three variables, say $x_1, x_2, x_3$, we may assume that $\{1,2\} \not\in G$ and $\{1,3\} \in G$. Then $x_1x_2 \in I_G$ and $f \in P_{13}^2$. Since $P_{13}$ does not contain $x_1,x_3$, $f$ is divisible by $x_2^2$. Hence $f$ is divisible by $x_1x_2^2x_3$. If $\{2,3\} \not\in G$, then $x_2x_3 \in I_G$, which implies $f \in I_G^2$. If $\{2,3\} \in G$, then $f$ is also divisible by $x_1^2$. Hence $f$ is divisible by $x_1^2x_2^2$, which implies $f \in I_G^2$.\par
If $f$ involves more than four variables, we suppose that $f$ is divisible by $x_1x_2x_3x_4$.
By the assumption on $G$ we may assume that $\{1,2\} \in G$ and $\{1,3\} \not\in G$. If $\{2,4\} \not\in G$, then $x_1x_3,x_2x_4 \in I_G$, hence $f \in I_G^2$. If $\{2,4\} \in G$, then $\{1,4\} \not\in G$. From this it follows that $\{3,4\} \in G$, hence we must have $\{2,3\} \not\in G$. Therefore, $x_1x_4, x_2x_3 \in I_G$, which implies $f\in I_G^2$.
\end{proof}

For $m \ge 3$, we shall see that the condition $I_G^{(m)} = I_G^m$ implies $n \le 4$ so that we get the following criterion.

\begin{Theorem} \label{equa 2}
Let $n \ge 4$. For $m \ge 3$, $I_G^{(m)} = I_G^m$ if and only if $n =4$ and $G$ is a path or a cycle or the union of two disjoint edges.
\end{Theorem}

\begin{proof}
Assume that $I_G^{(m)} = I_G^m$ for some $m \ge 3$. 
If $n \ge 5$, then $x_1x_2x_3x_4x_5 \in P_{ij}^3$ for all $i,j$. Hence $x_1x_2x_3x_4x_5 \in I_G^{(3)}$. We may assume that $\{1,2\} \not\in G$. Then $x_1x_2 \in I_G$. Therefore, 
$$x_1^{m-2}x_2^{m-2}x_3x_4x_5 = (x_1x_2)^{m-3}(x_1x_2x_3x_4x_5)\in I^{m-3}I^{(3)} \subseteq I_G^{(m)}.$$
This monomial has degree $2m-1$ so that it can not belong to $I_G^m$, a contradiction. Thus, $n = 4$. By Lemma \ref{triangle}, this implies that $G$ is a path or a cycle  or the union of two disjoint edges. \par
Conversely, if $G$ is one of the graphs in the assertion, then $I_G$ is the edge ideal of the complementary graph of $G$ which consists of the non-edges of $G$. This graph is a bipartite graph. Hence $I_G^{(m)} = I_G^m$ for all $m \ge 2$ by \cite[Theorem 5.9]{SiVV}.
\end{proof}

Now, we can give the following criteria for the Cohen-Macaulayness of $I_G^m$.

\begin{Corollary}
Let $n \ge 4$. Then $I_G^2$ is a Cohen-Macaulay ideal if and only if $G$ is a cycle of length 4,5.
\end{Corollary}

\begin{proof}
This follows from Theorem \ref{CM 1} and Theorem \ref{equa 1}.
\end{proof}

\begin{Corollary}
Let $n \ge 4$. For $m \ge 3$,  $I_G^m$ is a Cohen-Macaulay ideal if and only if $G$ is a cycle of length 4.
\end{Corollary}

\begin{proof}
This follows from Theorem \ref{CM 2} and Theorem \ref{equa 2}.
\end{proof}

If $n = 3$ or if $G$ is a cycle of length 4, $I_G$ is a complete intersection.
Therefore, $I_G^m$ is a Cohen-Macaulay ideal for some $m \ge 3$ if and only if $I_G$ is a complete intersection.
This fact does not hold for $m = 2$ because $I_G$ is not a complete intersection if $G$ is a cycle of length 5.
\smallskip

We can also use the above results to study the vertex cover algebras of certain simplicial complexes.
\smallskip 

Let $\Delta$ be a simplicial complex on the vertex set $\{1,...,n\}$. An integer vector $\c = (c_1, \ldots, c_n) \in
\NN^n$ is called a $m$-cover of $\Delta$ if $\sum_{i \in F} c_i \geq m$ for all facets $F$ of
$\Delta$. Let $A_m(\Delta)$ denote the $k$-vector space generated by all monomials $x^\c t^m$ such that $\c$ is a $m$-cover of $\Delta$, where $t$ is a new variable. Then
\[
A(\Delta) :=  \bigoplus_{m \geq 0}A_m(\Delta),
\]
is a graded $S$-algebra. We call $A(\Delta)$ the  vertex cover algebra of $\Delta$ \cite{HHT}.

Vertex cover algebra has an interesting algebraic interpretation. Let ${\mathcal F}(\Delta)$ denote the set of the facets of $\Delta$. Then 
$$I^*(\Delta)  :=\bigcap_{F\in {\mathcal F}(\Delta)}P_{\bar F},$$
where $\bar F$ denotes the complement of $F$.
Then $A(\Delta)$ is the symbolic Rees algebra of $I^*(\Delta)$. 
It is shown in \cite{HHT} that $A(\Delta)$  is a finitely generated, graded and normal 
Cohen-Macaulay $S$-algebra.
\smallskip

It is of great interest to know when $A(\Delta)$ is a standard graded algebra, that is, when $A(\Delta)$  is generated over $S$ by forms of degree 1. This is equivalent to the condition $I^*(\D)^{(m)} = I^*(\D)^m$ for all $m \ge 1$ \cite{HHT} and can be described in terms of the max-flow min-cut property in integer programming \cite{HHTZ}. Note that this condition is always satisfied if $\dim \D = 0$ ($I^*(\D)$ is a principal ideal).

\begin{Corollary} 
Let $\D$ be a pure $(n-3)$-dimensional simplicial complex on $n$ vertices, $n \ge 4$. Then $A(\D)$ is a standard graded algebra if and only if $n = 4$ and  $\D$ is a path or a cycle of length $4$ or  the union of two disjoint edges. 
\end{Corollary}

\begin{proof}
The assumption on $\D$ implies that the complements of the facets of $\D$ form a graph $G$.
Since $I^*(\D) = I_G$, applying Theorem \ref{equa 1} and Theorem \ref{equa 2} to $G$ we see that 
$A(\D)$ is a standard graded algebra if and only if $n = 4$ and $\D$ is a path or a cycle or  the union of two disjoint edges. 
\end{proof}

\end{document}